\documentclass[12pt]{amsart}
\usepackage{amssymb}
\usepackage[all]{xy}
\usepackage[toc,page,title,titletoc,header]{appendix}

\begin{document}
\theoremstyle{plain}
\newtheorem{thm}{Theorem}[section]
\newtheorem*{thm1}{Theorem 1}
\newtheorem*{thm1.1}{Theorem 1.1}
\newtheorem*{thmM}{Main Theorem}
\newtheorem*{thmA}{Theorem A}
\newtheorem*{thm2}{Theorem 2}
\newtheorem{lemma}[thm]{Lemma}
\newtheorem{lem}[thm]{Lemma}
\newtheorem{cor}[thm]{Corollary}
\newtheorem{pro}[thm]{Proposition}
\newtheorem{propose}[thm]{Proposition}
\newtheorem{variant}[thm]{Variant}
\theoremstyle{definition}
\newtheorem{notations}[thm]{Notations}
\newtheorem{rem}[thm]{Remark}
\newtheorem{rmk}[thm]{Remark}
\newtheorem{rmks}[thm]{Remarks}
\newtheorem{defi}[thm]{Definition}
\newtheorem{exe}[thm]{Example}
\newtheorem{claim}[thm]{Claim}
\newtheorem{ass}[thm]{Assumption}
\newtheorem{prodefi}[thm]{Proposition-Definition}
\newtheorem{que}[thm]{Question}
\newtheorem{con}[thm]{Conjecture}

\newtheorem*{dmlcon}{Dynamical Mordell-Lang Conjecture}
\newtheorem*{condml}{Dynamical Mordell-Lang Conjecture}
\numberwithin{equation}{section}
\newcounter{elno}                % This to number lists
\def\points{\list
{\hss\llap{\upshape{(\roman{elno})}}}{\usecounter{elno}}}
\let\endpoints=\endlist
\newcommand{\SH}{\rm SH}
\newcommand{\Tan}{\rm Tan}
\newcommand{\res}{\rm res}
\newcommand{\Om}{\Omega}
\newcommand{\om}{\omega}
\newcommand{\la}{\lambda}
\newcommand{\mc}{\mathcal}
\newcommand{\mb}{\mathbb}
\newcommand{\surj}{\twoheadrightarrow}
\newcommand{\inj}{\hookrightarrow}
\newcommand{\zar}{{\rm zar}}
\newcommand{\Exc}{\rm Exc}
\newcommand{\an}{{\rm an}}
\newcommand{\red}{{\rm red}}
\newcommand{\codim}{{\rm codim}}
\newcommand{\Supp}{{\rm Supp}}
\newcommand{\rank}{{\rm rank}}
\newcommand{\Ker}{{\rm Ker \ }}
\newcommand{\Pic}{{\rm Pic}}
\newcommand{\Div}{{\rm Div}}
\newcommand{\Hom}{{\rm Hom}}
\newcommand{\im}{{\rm im}}
\newcommand{\Spec}{{\rm Spec \,}}
\newcommand{\Nef}{{\rm Nef \,}}
\newcommand{\Frac}{{\rm Frac \,}}
\newcommand{\Sing}{{\rm Sing}}
\newcommand{\sing}{{\rm sing}}
\newcommand{\reg}{{\rm reg}}
\newcommand{\Char}{{\rm char}}
\newcommand{\Tr}{{\rm Tr}}
\newcommand{\ord}{{\rm ord}}
\newcommand{\id}{{\rm id}}
\newcommand{\NE}{{\rm NE}}
\newcommand{\Gal}{{\rm Gal}}
\newcommand{\Min}{{\rm Min \ }}
\newcommand{\Max}{{\rm Max \ }}
\newcommand{\Alb}{{\rm Alb}\,}
\newcommand{\GL}{{\rm GL}\,}        % For the general linear group
\newcommand{\PGL}{{\rm PGL}\,}
\newcommand{\Bir}{{\rm Bir}}
\newcommand{\Aut}{{\rm Aut}}
\newcommand{\End}{{\rm End}}
\newcommand{\Per}{{\rm Per}\,}
\newcommand{\ie}{{\it i.e.\/},\ }
\newcommand{\niso}{\not\cong}
\newcommand{\nin}{\not\in}
\newcommand{\soplus}[1]{\stackrel{#1}{\oplus}}
\newcommand{\by}[1]{\stackrel{#1}{\rightarrow}}
\newcommand{\longby}[1]{\stackrel{#1}{\longrightarrow}}
\newcommand{\vlongby}[1]{\stackrel{#1}{\mbox{\large{$\longrightarrow$}}}}
\newcommand{\ldownarrow}{\mbox{\Large{\Large{$\downarrow$}}}}
\newcommand{\lsearrow}{\mbox{\Large{$\searrow$}}}
\renewcommand{\d}{\stackrel{\mbox{\scriptsize{$\bullet$}}}{}}
\newcommand{\dlog}{{\rm dlog}\,}    % For dlog
\newcommand{\longto}{\longrightarrow}
\newcommand{\vlongto}{\mbox{{\Large{$\longto$}}}}
\newcommand{\limdir}[1]{{\displaystyle{\mathop{\rm lim}_{\buildrel\longrightarrow\over{#1}}}}\,}
\newcommand{\liminv}[1]{{\displaystyle{\mathop{\rm lim}_{\buildrel\longleftarrow\over{#1}}}}\,}
\newcommand{\norm}[1]{\mbox{$\parallel{#1}\parallel$}}
\newcommand{\boxtensor}{{\Box\kern-9.03pt\raise1.42pt\hbox{$\times$}}}
\newcommand{\into}{\hookrightarrow}
\newcommand{\image}{{\rm image}\,}
\newcommand{\Lie}{{\rm Lie}\,}      % For Lie algebra of groups
\newcommand{\CM}{\rm CM}
\newcommand{\sext}{\mbox{${\mathcal E}xt\,$}}  % For sheaf Ext
\newcommand{\shom}{\mbox{${\mathcal H}om\,$}}  %For sheaf Hom
\newcommand{\coker}{{\rm coker}\,}  % For the cokernel of a morphism
\newcommand{\sm}{{\rm sm}}
\newcommand{\pgcd}{\text{pgcd}}
\newcommand{\trd}{\text{tr.d.}}
\newcommand{\tensor}{\otimes}
\renewcommand{\iff}{\mbox{ $\Longleftrightarrow$ }}
\newcommand{\supp}{{\rm supp}\,}
\newcommand{\ext}[1]{\stackrel{#1}{\wedge}}
\newcommand{\onto}{\mbox{$\,\>>>\hspace{-.5cm}\to\hspace{.15cm}$}}
\newcommand{\propsubset}
{\mbox{$\textstyle{
\subseteq_{\kern-5pt\raise-1pt\hbox{\mbox{\tiny{$/$}}}}}$}}
% Skriptbuchstaben
\newcommand{\sA}{{\mathcal A}}
\newcommand{\sB}{{\mathcal B}}
\newcommand{\sC}{{\mathcal C}}
\newcommand{\sD}{{\mathcal D}}
\newcommand{\sE}{{\mathcal E}}
\newcommand{\sF}{{\mathcal F}}
\newcommand{\sG}{{\mathcal G}}
\newcommand{\sH}{{\mathcal H}}
\newcommand{\sI}{{\mathcal I}}
\newcommand{\sJ}{{\mathcal J}}
\newcommand{\sK}{{\mathcal K}}
\newcommand{\sL}{{\mathcal L}}
\newcommand{\sM}{{\mathcal M}}
\newcommand{\sN}{{\mathcal N}}
\newcommand{\sO}{{\mathcal O}}
\newcommand{\sP}{{\mathcal P}}
\newcommand{\sQ}{{\mathcal Q}}
\newcommand{\sR}{{\mathcal R}}
\newcommand{\sS}{{\mathcal S}}
\newcommand{\sT}{{\mathcal T}}
\newcommand{\sU}{{\mathcal U}}
\newcommand{\sV}{{\mathcal V}}
\newcommand{\sW}{{\mathcal W}}
\newcommand{\sX}{{\mathcal X}}
\newcommand{\sY}{{\mathcal Y}}
\newcommand{\sZ}{{\mathcal Z}}
% Sonderbuchstaben mit Doppellinie
\newcommand{\A}{{\mathbb A}}
\newcommand{\B}{{\mathbb B}}
\newcommand{\C}{{\mathbb C}}
\newcommand{\D}{{\mathbb D}}
\newcommand{\E}{{\mathbb E}}
\newcommand{\F}{{\mathbb F}}
\newcommand{\G}{{\mathbb G}}
\newcommand{\HH}{{\mathbb H}}
\newcommand{\I}{{\mathbb I}}
\newcommand{\J}{{\mathbb J}}
\newcommand{\M}{{\mathbb M}}
\newcommand{\N}{{\mathbb N}}
\renewcommand{\P}{{\mathbb P}}
\newcommand{\Q}{{\mathbb Q}}
\newcommand{\R}{{\mathbb R}}
\newcommand{\T}{{\mathbb T}}
\newcommand{\U}{{\mathbb U}}
\newcommand{\V}{{\mathbb V}}
\newcommand{\W}{{\mathbb W}}
\newcommand{\X}{{\mathbb X}}
\newcommand{\Y}{{\mathbb Y}}
\newcommand{\Z}{{\mathbb Z}}

\newcommand{\fix}{\mathrm{Fix}}

%%%%%%%%%%%%%%%%%%%%%%%%%%%%%%%%%%%%%%%%%%%%%%%%%%%%%%%%%%%%%%
%%%%%%%%%%%%%%%%%%%%%%%%%%%%%%%%%%%%%%%%%%%%%%%%%%%%%%%%%%%%%%
\title[]{The existence of Zariski dense orbits for polynomial endomorphisms of the affine plane}

\author{Junyi Xie}

\address{Universit\'e de Rennes I
  Campus de Beaulieu,
  b\^atiment 22-23,
  35042 Rennes cedex
  France}

\email{junyi.xie@univ-rennes1.fr}

%\date{\today}

\bibliographystyle{plain}

\maketitle

\begin{abstract}
In this paper we prove the following theorem.
Let $f$ be a dominant polynomial endomorphism of the affine plane defined over an algebraically closed field of characteristic $0$. If there is no nonconstant invariant rational function under $f$,
then there exists a closed point in the plane whose orbit under $f$ is Zariski dense.

This result gives us a positive answer to a conjecture proposed by Medvedev and Scanlon, by Amerik, Bogomolov and Rovinsky, and by Zhang, for polynomial endomorphisms of the affine plane.
\end{abstract}

%
%\tableofcontents

\renewcommand{\thefootnote}{}
\footnotetext{The author is supported by the labex CIMI.}

\section{Introduction}
Denote by $k$ an algebraically closed field of characteristic $0$.

The aim of this paper is to prove
\begin{thm}\label{thmconzhangatwo}Let $f:\mathbb{A}^2_k\rightarrow \mathbb{A}^2_k$ be a dominant polynomial endomorphism.
If there are no nonconstant
rational functions $g$ satisfying $g\circ f= g$,
there exists a point $p\in \mathbb{A}^2(k)$ such that the orbit $\{f^n(p)|\,\,n\geq 0\}$ of $p$ is Zariski dense in $\mathbb{A}^2_k$.
\end{thm}
We cannot ask $g$ in Theorem \ref{thmconzhangatwo} to be a polynomial.
Indeed, let $P(x,y)$ be a polynomial which is neither zero nor a root of unity. Let $f:\A^2_k\to \A^2_k$ be the endomorphism defined by $(x,y)\mapsto (P(x,y)x,P(x,y)y)$. It is easy to see that $g\circ f=g$ if $g=y/x$, but there does not exist any polynomial $h$ satisfying $h\circ f=h$.

\medskip

The following conjecture was proposed by Medvedev and Scanlon \cite[Conjecture 5.10]{Medvdevv1} and also by Amerik, Bogomolov and Rovinsky \cite{E.Amerik2011}
\begin{con}\label{conexistszdo}Let $X$ be a quasi-projective variety over $k$ and $f:X\rightarrow X$ be a dominant endomorphism for which there exists no nonconstant
rational function $g$ satisfying $g\circ f=g$.
Then there exists a point $p\in X(k)$ whose orbit is Zariski dense in $X$.
\end{con}
Conjecture \ref{conexistszdo} strengthens the following conjecture of Zhang \cite{zhang}.
\begin{con}\label{conexistszdozhang}Let $X$ be a projective variety and $f:X\rightarrow X$ be an endomorphism defined over $k$. If there exists an ample line bundle $L$ on $X$ satisfying $f^*L=L^{\otimes d}$ for some integer $d>1$,
then there exists a point $p\in X(k)$ whose orbit $\{f^n(p)|\,\,n\geq 0\}$ is Zariski dense in $X(k)$.
\end{con}
Theorem \ref{thmconzhangatwo} setlles Conjecture \ref{conexistszdo} for polynomial endomorphisms of $\mathbb{A}^2_{k}$.

When $k$ is uncountable, Conjecture \ref{conexistszdo} was proved by Amerik and Campana \cite{Amerik2008}.
In \cite{Fakhruddin2014}, Fakhruddin proved Conjecture \ref{conexistszdo} for  $generic^1$\footnote{An endomorphism $f:\mathbb{P}^N_{k}\rightarrow \mathbb{P}^N_{k}$ satisfying $f^*O_{\mathbb{P}^N_{k}}(1)=O_{\mathbb{P}^N_{k}}(d)$ is said to be generic if it conjugates by a suitable linear automorphism on $\P^N_{k}$ to an endomorphism $[x_0:\cdots:x_N]\mapsto [\sum_{|I|=d}a_{0,I}x^I:\cdots: \sum_{|I|=d}a_{N,I}x^I]$ where the set $\{a_{i,I}\}_{0\leq i\leq N, |I|=d}$ is algebraically independent over $\bar{\mathbb{Q}}$.} endomorphisms on projective spaces over arbitrary algebraically closed fields $k$ of characteristic zero.
In \cite{Xie2015}, the author proved Conjecture \ref{conexistszdo} for birational surface endomorphisms with dynamical degree great than $1$.  Recently in \cite{Medvdev},  Medvedev and Scanlon proved Conjecture \ref{conexistszdo} when $f:=(f_1(x_1),\cdots,f_N(x_N))$ is an endomorphism of $\mathbb{A}^N_{k}$
where the $f_i$'s are one-variable polynomials defined over $k$.

We mention that in \cite{Amerik}, Amerik proved that there exists a nonpreperiodic algebraic point when  $f$ is of infinite order. In \cite{Bell}, Bell, Ghioca and Tucker proved that if $f$ is an automorphism, then there exists a
subvariety of codimension $2$ whose orbit under $f$ is Zariski dense.

We note that Conjecture 1.2 is not true in the case when $k$ is the algebraic
closure of a finite field, since in this case all orbits of $k$-points are finite.

%However,
%I suspect that when $k$ is transcendental over a finite field, Conjecture 1.2 should
%hold i.e.
%
%\begin{con}
%Let $k$ be an algebraically closed field of characteristic $p>0$
%which is a transcendental extension over $\F_p$. Let $X$ be a quasi-projective variety
%over $k$ and $f: X \to X$ be a dominant endomorphism for which there exists no
%nonconstant rational function $g$ satisfying $g\circ f=g$. Then there exists a point
%$x\in X(k)$ whose orbit is Zariski dense in $X$.
%\end{con}
%We will study this conjecture in the forthcoming work \cite{Xie2016charp}.

\medskip

Our proof of Theorem \ref{thmconzhangatwo} is based on the valuative techniques developed in \cite{Favre2004,Favre2007,Favre2011,Xiec}. Here is an outline of the proof.

For simplicity, suppose that $f$ is a dominant polynomial map $f:=(F(x,y),G(x,y))$ defined over $\Z$.

 When $f$ is birational, our Theorem \ref{thmconzhangatwo} is essentially proved in \cite{Xie2015}. So we may suppose that $f$ is not birational.

 By \cite{Favre2011}, there exists a projective compactification $X$ of $\A^2$ for which the induced
map by $f$ at infinity is \emph{algebraically stable} i.e. it does not contract any
curve to a point of indeterminacy. Moreover, we can construct an "attracting locus" at infinity,  in the sense that
\begin{points}
\item
either there exists a superattracting fixed point $q\in X\setminus \A^2$ such that there is no branch of curve at $q$ which is periodic under $f$;
\item or there exists an irreducible component $E\in X\setminus \A^2$ such that $f^*E=dE+F$ where $d\geq 2$ and $F$ is an effective divisor supported by $X\setminus \A^2.$
\end{points}

In Case (i), we can find a point $p\in \A^2(\bar{\Q})$, near $q$ w.r.t. the euclidean topology. Then we have $\lim_{n\to \infty} f^n(p)=q$. It is easy to show that the orbit of $p$ is Zariski dense in $\A^2$.

In Case (ii), $E$ is defined over $\Q.$ There exists a prime number ${\mathfrak{p}}\geq 3$, such that $f_{\mathfrak{p}}|_{E_{\mathfrak{p}}}$ is dominant where $f_{\mathfrak{p}}:= f\,  \rm{mod }\,  {\mathfrak{p}}$ and $E_{\mathfrak{p}}:=E \, \rm{mod }\,  {\mathfrak{p}}$.

We first treat the case $f^n|_E\neq \id$ for all $n\geq 1$.
 After replacing $f$ by a suitable iterate, we may find a fixed point $x\in E_{\mathfrak{p}}$ such that $df_{\mathfrak{p}}(x)=1$ in $\bar{\F}_{\mathfrak{p}}.$ Denote by $U$ the ${\mathfrak{p}}$-adic open set of $X(\Q_{\mathfrak{p}})$ consisting the points $y$ such that $y \, \rm{mod }\,  {\mathfrak{p}}=x$. Then $U$ is fixed by $f$. By \cite[Theorem 1]{Poonen2014}, all the preperiodic points in $U\cap E$ are fixed. Moreover $\bigcap_{n\geq 0}f^n(U)=U\cap E.$ Denote by $S$ the set of fixed points in $U\cap E$. Then $S$ is finite.
If $S$ is empty, pick a point $p\in \A^2(\bar{\Q}\cap \Q_p)\cap U$, it is easy to see that the orbit of $p$ is Zariski dense in $\A^2$.
If $S$ is not empty, by \cite[Theorem 3.1.4]{Abate2001}, at each point $q_i\in S$, there exists at most one algebraic curve $C_i$ passing through $q_i$ which is preperiodic. Set $C_i=\emptyset$ if no such curve does exists. We have that $C_i$ is fixed. Pick a point $p\in \A^2(\bar{\Q}\cap \Q_p)\setminus C_1$ very closed to $q_1$. We can show that the orbit of $p$ is Zariski dense in $\A^2$.

Next, we treat the case $f|_E=\id$. By \cite[Theorem 3.1.4]{Abate2001}, at each point $q\in E$, there exists at most one algebraic curve $C_q$ passing through $q$ which is preperiodic. Set $C_q=\emptyset$, if such curve does not exists. We have that $C_q$ is fixed and transverse to $E$. If $C_q=\emptyset$ for all but finitely many $q\in E$, there exists $q\in E$ and a $\mathfrak{p}$-adic neighborhood $U$ of $q$ such that for any point $y\in U\cap E$, $C_y=\emptyset$, $f(U)\subseteq U$ and $\bigcap_{\infty}f^n(U)=U\cap E.$ Then for any point $p\in \A^2(\bar{\Q}\cap \Q_p)\cap U$, the orbit of $p$ is Zariski dense in $\A^2$. Otherwise there exists a sequence of points $q_i\in E$, such that $C_i:=C_{q_i}$ is an irreducible curve. Since $f|_{C_i}$ is an endomorphism of $C_i\cap \A^2$ of degree at least two, $C_i$ has at most two branches at infinity. Since $C_{i}$ is transverse to $E$ at $q_i$, we can bound the intersection number $(E\cdot C_i)$. By the technique developed in \cite{Xiec}, we can also bound the intersection of $C_i$ with the other irreducible components of $X\setminus \A^2.$ Then we bound the degree of $C_i$ which allows us to construct a nonconstant invariant rational function.

\medskip

The article is organized in 2 parts.

In Part \ref{part1}, we gather some results on
the geometry and dynamics at infinity and metrics on projective varieties defined
over a valued field. We first introduce the valuative tree at
infinity in Section \ref{sectiontvaluationtreeinfty}, and then we recall the main properties of the action of a
polynomial map on the valuation space in Section \ref{sectionbgdopm}. Next we introduce the Green function on the valuative tree for a polynomial endomorphism in Section \ref{sectiongreen}.
Finally we give background information
on metrics on projective varieties defined
over a valued field in Section \ref{secmop}.

In Part \ref{part2}, we prove Theorem \ref{thmconzhangatwo}. We first prove it in some special cases in Section \ref{sectionatt}. In most of these cases, we find a Zariski dense orbit in some attracting locus. Then we study totally invariant curves in Section \ref{sectiontot} and prove Theorem \ref{thmconzhangatwo} when there are infinitely many such curves. Finally we finish the proof of Theorem \ref{thmconzhangatwo} in Section \ref{sectionproof}.

\section*{Acknowledgement}
Part of this work was done during the PhD study of the author with Charles
Favre in Ecole polytechnique. I would like to thank Charles Favre for his constant
support during this time. The referee carefully read a first version of this paper and his/her comments
greatly improved it. I also thank Serge Cantat and St\'ephane Lamy for 
useful discussions. I thank Vladimir Popov and Dragos Ghioca for their comments
on the first version of this paper. I finally thank Matteo Ruggiero for answering
my question about the theorem of Hadamard-Perron \cite[Theorem 3.1.4]{Abate2001}.

%
%
%Part of this work was done during the PhD study of the author with Charles Favre in Ecole polytechnique. I would like to thank Charles Favre for his constant support during this time. I also thank Serge Cantat and St\'ephane Lamy for those useful discussions. I finally thank Matteo Ruggiero for answering my question about the theorem of Hadamard-Perron \cite[Theorem 3.1.4]{Abate2001}.

\part{Preliminaries}\label{part1}

In this part, we denote by $k$ an algebraically closed field of characteristic zero. We also fix affine coordinates on $\A^2_k = \Spec k [x,y]$.

\section{The valuative tree at infinity}\label{sectiontvaluationtreeinfty}
 We refer to \cite{Jonsson} for details, see also \cite{Favre2004,Favre2007,Favre2011}.
\subsection{The valuative tree at infinity}

In this article by a valuation on a unitary $k$-algebra $R$ we shall understand a function
$v : R \to \mathbb{R}\cup \{+\infty\}$ such that the restriction of $v$ to $k^* = k - \{ 0 \}$ is
constant equal to $0$, and $v$ satisfies $v(fg) = v(f) + v(g)$ and $v(f +g) \ge \min \{v(f), v(g) \}$.
It is usually referred to as a semivaluation in the literature, see \cite{Favre2004}. We will however make a slight abuse of notation and call it a valuation.

\smallskip

Denote by $V_{\infty}$ \index{$V_{\infty}$}the space of all normalized valuations centered at infinity \index{normalized valuation centered at infinity} i.e. the set of valuations $v:k[x,y]\rightarrow \mathbb{R}\cup \{+\infty\}$ satisfying  $\min\{v(x),v(y)\}=-1$. The topology on $V_{\infty}$ is defined to be the weakest topology making the map $v\mapsto v(P)$ continuous for every $P\in k[x,y]$.

\smallskip

The set $V_{\infty}$ is equipped with a {\em partial ordering} defined by $v\leq w$ if and only if
$v(P)\leq w(P)$ for all $P\in k[x,y]$. Then  $-\deg: P \mapsto -\deg (P)$ is the unique minimal element.

\smallskip

Given any valuation $v\in V_\infty\setminus \{-\deg\}$, the set
$ \{ w \in V_\infty, \,| - \deg \le w \le v \}$ is isomorphic as a poset to the real segment
$[0,1]$ endowed with the standard ordering. In other words, $(V_\infty, \le)$ is a rooted tree in the sense of \cite{Favre2004,Jonsson}.

Given any two valuations $v_1,v_2\in V_{\infty}$,
there is a unique valuation in $V_{\infty}$ which is maximal in the set $\{v\in V_{\infty}|\,\,v\leq v_1 \text{ and } v\leq v_2\}.$ We denote it by  $v_1\wedge v_2$.

The segment $[v_1, v_2]$ is by definition the union of $\{w \,| \, v_1\wedge v_2 \le w \le v_1\}$
and $\{w \,| \, v_1\wedge v_2 \le w \le v_2\}$.

\smallskip

Pick any valuation  $v\in V_\infty$. We say that two points $v_1, v_2$
lie in the same direction at $v$ if the segment $[v_1, v_2]$ does not contain $v$.
A {\em direction} (or a tangent vector) at $v$ is an equivalence class for this relation.
We write $\Tan_v$ for the set of directions at $v$.

\smallskip

When $\Tan_v$ is a singleton, then $v$ is called an endpoint\index{endpoint}. In $V_\infty$, the set of endpoints is exactly the set of all maximal valuations. %This set is dense in $V_{\infty}.$
When $\Tan_v$ contains exactly two directions, then $v$ is said to be regular.
When $\Tan_v$ has more than three directions, then $v$ is a branch point.

\smallskip

Pick any $v\in V_\infty$.
For any tangent vector $\vec{v}\in \Tan_v$, denote by $U(\vec{v})$ the subset of those elements in $V_\infty$ that determine $\vec{v}$. This is an open set whose boundary is reduced to the singleton $\{v\}$. If $v\neq -\deg$, the complement of $\{w\in V_\infty, \,| w \ge v\}$ is equal to $U(\vec{v}_0)$ where $\vec{v}_0$ is the tangent vector determined by $-\deg$.

It is a fact that finite intersections of open sets of the form $U(\vec{v})$ form a basis for the topology of $V_\infty$.

%
%\smallskip
%
%The \emph{convex hull} of  any subset $S\subset V_\infty$ is defined as the set of valuations $v\in V_\infty$ such that
%there exists a pair $v_1 , v_2 \in S$ with $v\in [v_1, v_2]$.
%
%A \emph{finite subtree} of $V_\infty$ is, by definition, the convex hull of a finite collection of points in $V_\infty$. A point in a finite subtree  $T\subseteq V_{\infty}$ is said to be an end point if it is extremal in $T.$
%

\subsection{Compactifications of $\mathbb{A}^2_k$}

A {\em compactification} of $\A^2_k$ is the data of a projective surface $X$ together with an open immersion $\A^2_k \to X$ with dense image.

A compactification $X$ dominates another one $X'$ if the canonical birational map $X \dashrightarrow X'$ induced by the inclusion of $\A^2_k$ in both surfaces is in fact a regular map.

The category $\mathcal{C}$ of all compactifications of $\A^2_k$ forms an inductive system for the relation of domination.

\smallskip

Recall that we have fixed affine coordinates on $\A^2_k = \Spec k[x,y]$.
We let $\P^2_k$ be the standard compactification of $\A^2_k$ and denote by
$l_\infty := \P^2_k \setminus \A^2_k$ the line at infinity in the projective plane.

An {\em admissible compactification} of $\A^2_k$ is by definition a smooth projective surface $X$ endowed with a birational morphism $\pi_X : X \to \P^2_k$ such that $\pi_X$ is an isomorphism over $\A^2_k$ with the embedding $\pi^{-1}|_{\A_k^2}:\A_k^2\to X$. Recall that $\pi_X$ can then be decomposed as a finite sequence of point blow-ups.

We shall denote by $\mathcal{C}_0$ the category of all admissible compactifications. It is also an inductive system for the relation of domination. Moreover $\mathcal{C}_0$ is a subcategory of $\mathcal{C}$ and for any compactification $X\in \mathcal{C}$, there exists that $X'\in \mathcal{C}_0$ dominates $X$.

\subsection{Divisorial valuations}
Let $X\in \mathcal{C}$ be a compactification of $\A^2_k=\Spec k[x,y]$ and $E$ be an irreducible component of $X\setminus \mathbb{A}^2$. Set $b_E:=-\min\{\ord_E(x),\ord_E(y)\}$  and $v_E:=b_E^{-1}\ord_E$. Then we have $v_E\in V_{\infty}$.

By Poincar\'e Duality there exists a unique
{\em dual divisor} \index{dual divisor} $\check{E}$  of $E$ defined as the unique divisor supported on $X\setminus \A^2$ such that $(\check{E}\cdot F)=\delta_{E,F}$ for all irreducible components $F$ of $X\setminus \A^2$.

\begin{rem}\label{rembegeo}
Recall that $l_{\infty}$ is the line at infinity of $\mathbb{P}^2$. Let $s$ be a formal
curve centered at some point $q\in l_{\infty}.$ Suppose that the strict transform of $s$ in $X$ intersects $E$ transversally at a point in $E$ which is smooth in $X\setminus \A^2.$ Then we have $(s\cdot l_{\infty})=b_E$.
\end{rem}

%XXX define the dual divisor XXX

\subsection{Classification of valuations}
There are four kinds of valuations in $V_{\infty}$. The first one corresponds to the {\em divisorial valuations} which we have defined above. We now describe the three remaining types of valuations.
\subsubsection*{Irrational valuations}
Consider any two irreducible components $E$ and $E'$ of $X\setminus \A^2_k$ for some compactification $X\in \mathcal{C}$ of $\mathbb{A}^2_k$
intersecting at a point $p$. There exists local coordinates $(z,w)$ at $p$ such that $E=\{z=0\}$
and $E'=\{w=0\}$. To any pair $(s,t)\in (\mathbb{R}^{+})^{2}$ satisfying  $sb_E+tb_{E'}=1$, we attach the valuation $v$ defined on the
ring $O_p$ of germs at $p$ by
$v(\sum a_{ij}z^iw^j)=\min \{si+tj|\,\, a_{ij}\neq 0\}.$
Observe that it does not
depend on the choice of coordinates. By first extending $v$ to the common fraction
field $k(x,y)$ of $O_p$ and $k[x,y]$, then restricting it to $k[x,y]$, we obtain a valuation in $V_{\infty}$, called
{\em quasimonomial}\index{quasimonomial}.  It is divisorial if and only if either $t=0$ or the ratio $s/t$ is a
rational number. Any divisorial valuation is quasimonomial. An {\em irrational valuation} is
by definition a nondivisorial quasimonomial valuation.
\subsubsection*{Curve valuations}
Recall that $l_{\infty}$ is the line at infinity of $\mathbb{P}^2_k$.
For any formal curve $s$ centered at some point $q\in l_{\infty}$, denote by $v_s$\index{$v_s$} the valuation defined by $P\mapsto (s\cdot l_{\infty})^{-1}\ord_{\infty}(P|_s)$. Then we have $v_s\in V_{\infty}$ and call it a {\em curve valuation}\index{curve valuation}.

\medskip

Let $C$ be an irreducible curve in $\P^2_k$. For any point $q\in C\cap l_{\infty}$, denote by $O_{q}$ the local ring at $q$, $m_q$ the maximal ideal of $O_q$ and $I_C$ the ideal of height 1 in $O_q$ defined by $C$. Denote by $\widehat{O}_{q}$ the completion of $O_{q}$ w.r.t. $m_q$, $\widehat{m}_q$ the completion of $m_q$ and $\widehat{I}_C$ the completion of $I_C$. For any prime ideal $\widehat{p}$ of height 1 containing $\widehat{I}_C$, the morphism $\Spec \widehat{O}_q/\widehat{p}\rightarrow\Spec \widehat{O}_{q}$ defines a formal curve centered at $q$. Such a formal curve is called a \emph{branch of $C$ at infinity}.

\subsubsection*{Infinitely singular valuations}
Let $h$ be a formal series of the form $h(z)=\sum_{k=0}^{\infty} a_kz^{\beta_k}$ with $a_k\in k^*$ and $\{\beta\}_k$ an increasing sequence of rational numbers
with unbounded denominators. Then $P\mapsto -\min\{\ord_{\infty}(x),\ord_{\infty}(h(x^{-1}))\}^{-1}\ord_{\infty}P(x,h(x^{-1}))$ defines
a valuation in $V_{\infty}$ called an infinitely singular valuation\index{infinitely singular valuation}.

%XXX define other kind of valuations (in particular irrational valutions) and say what are their position in the valuative tree XXX
%
%XXX put the definition of curve valuation here XXX
%
%If $C$ is a formal curve centered at infinity, denote by $v_C$ the valuation defined by $P\mapsto (C\cdot l_{\infty})\ord_{\infty}(P|_C)$ where $l_{\infty}$ is the line at infinity of $\mathbb{P}^2_k$. Then we have $v_C\in V_{\infty}$.
%
%
%There are four types of valuations in $V_{\infty}$: divisorial valuations, irational valuations, curve valuations and infinitely singular valuations.

\medskip

A valuation $v\in V_{\infty}$ is a branch point in $V_{\infty}$ if and only if it is diviorial, it is a regular point in $V_{\infty}$ if and only if it is an irrational valuation, and it is an endpoint in $V_{\infty}$ if and only if it is a curve valuation or an infinitely singular valuation.
Moreover, given any smooth projective compactification $X$ in which $v=v_E$, one proves that the
map sending an element $V_\infty$ to its center in $X$ induces a  map $\Tan_v \to E$ that is a bijection.

\subsection{Parameterizations}
The {\em skewness} function\index{skewness function}
$\alpha:V_{\infty}\rightarrow [-\infty,1]$ \index{$\alpha$} is the unique  function on $V_{\infty}$
that is continuous on segments, and satisfies
$$\alpha(v_E)=\frac1{b_E^2}(\check{E}\cdot \check{E})$$ where $E$ is any irreducible component of $X\setminus \mathbb{A}^2_k$ of any compactification $X$ of $\mathbb{A}^2_k$ and $\check{E}$ is the dual divisor of $E$ as defined above.

The skewness function is strictly decreasing, and upper semicontinuous.
%Therefor it  induces a metric $d_{V_\infty}$ on $V_\infty$ by setting
%$$d_{V_\infty}(v_1,v_2):=2\alpha(v_1\wedge v_2)-\alpha(v_1)-\alpha(v_2)$$
%for all $v_1, v_2\in V_\infty.$
%In particular, any segment in $V_\infty$ carries a canonical metric for which it becomes isometric to a real segment.
%
%\medskip
%
In an analogous way, one defines the {\em thinness} function\index{thinness function} $A:V_{\infty}\rightarrow [-2,\infty]$\index{$A:V_{\infty}\rightarrow [-2,\infty]$} as the unique, increasing, lower semicontinuous function on $V_{\infty}$ such that for any irreducible exceptional divisor $E$ in some  compactification $X\in \mathcal{C}$, we have $$A(v_E)=\frac1{b_E}\, (1 + \ord_E(dx\wedge dy))~.$$ Here we extend the differential form $dx\wedge dy$ to a rational differential form on $X$.

%%We refer to \cite[Section A.1]{Favre2007}.
%
%These parameterizations behave in the following way:
%\begin{points}
%\item when $v$ is a divisorial valuation, then  $\alpha(v)$ and $A(v)$ are rational numbers;
%\item when $v$ is an irrational valuation, then $\alpha(v), A(v)\in \mathbb{R}\setminus \mathbb{Q}$;
%\item when $v$ is a curve valuation, then  $\alpha(v)=-\infty$, and $A(v)=+\infty$;
%\item when $v$ is an infinitely singular valuation, then $\alpha(v)$ and $A(v)$ can be either finite or infinite.
%\end{points}

\subsection{Computation of local intersection numbers of curves at infinity}\label{sectionlocalinter}
Let $s_1$, $s_2$ be two different formal curves at infinity.
We denote by $(s_1\cdot s_2)$ the intersection number of these two formal curves in $\P^2$. This intersection number is always nonnegative, and it is positive if and only if $s_1$ and $s_2$ are centered at the same point.

Denote by $l_{\infty}$ the line at infinity in $\mathbb{P}^2_k.$ Denote by $v_{s_1},v_{s_2}$ the curve valuations associated to $s_1$ and $s_2.$

By \cite[Proposition 2.2]{Xiec}, we have $$(s_1\cdot s_2)=(s_1\cdot l_{\infty})(s_2\cdot l_{\infty})(1-\alpha(v_{s_1}\wedge v_{s_2})).$$

\section{Background on dynamics of polynomial maps}\label{sectionbgdopm}

Recall that the affine coordinates have been fixed, $\A^2_k = \Spec k[x,y]$.

\subsection{Dynamical invariants of polynomial mappings}
The (algebraic) degree of  a  dominant polynomial endomorphism   $f=(F(x,y),G(x,y))$ defined on $\mathbb{A}^2_k$ is  defined by
$$\deg(f):=\max\{\deg(F),\deg(G)\}~.$$\index{$\deg(f)$}
It is not difficult to show that the sequence $\deg(f^n)$ is sub-multiplicative, so that the limit
$\la_1(f):=\lim_{n\rightarrow \infty}(\deg(f^n))^{\frac{1}{n}}$\index{$\la_1(f)$} exists.
It is referred to as the {\em dynamical degree} of $f$, and it is a theorem of Favre and Jonsson that $\la_1(f)$ is always a quadratic integer, see \cite{Favre2011}.

\smallskip

The (topological) degree  $\la_2(f)$\index{$\la_2(f)$} of $f$ is defined to be the number of preimages of a general closed point in $\A^2(k)$; one has $\la_2(fg)=\la_2(f)\la_2(g).$

It follows from B\'ezout's theorem that $\la_2(f) \le \deg(f)^2$ hence
\begin{equation}\label{eq:ineq}
\la_1(f)^2\geq \la_2(f)~.
\end{equation}

The resonant case $\la_1(f)^2 = \la_2(f)$ is quite special and the following structure theorem for these maps is proven in \cite{Favre2011}.

\begin{thm}
Any polynomial endomorphism $f$ of $\A^2_k$ such that $\la_1(f)^2 = \la_2(f)$ is proper\footnote{We say a polynomial endomorphism $f$ of $\A^2_k$ is proper if it is a proper morphism between schemes. When $k=\C$, it means that the preimage of any compact set of $\C^2$ is compact.}, and we are in exactly one of the following two exclusive cases.
\begin{enumerate}
\item
$\deg(f^n) \asymp \la_1(f)^n$; there exists a compactification $X$ of $\A^2_k$ to which $f$ extends as a regular map $f : X \to X$.
\item
$\deg(f^n) \asymp n\la_1(f)^n$; there exist affine coordinates $x,y$
in which $f$ takes the form
$$f(x,y) = (x^l + a_1 x^{l-1}+ \ldots + a_l , A_0(x) y^l + \ldots + A_l(x))
$$
where $a_i \in k$ and $A_i \in k[x]$ with $\deg A_0 \ge1$, and $l = \la_1(f)$.
 \end{enumerate}
\end{thm}
\rem Regular endomorphisms as in (1) have been classified in \cite{Favre2011}.
\endrem
\subsection{Valuative dynamics}\label{subsectionvaldy}%\cite{Favre2011}
Any dominant polynomial endomorphism $f$ as in the previous section induces a natural map on the space of valuations at infinity in the following way.

\smallskip

For any $v\in V_\infty$ we set
$$d(f,v):=-\min\{v(F),v(G),0\}\geq 0~.$$
In this way, we get a non-negative continuous  decreasing function on $V_\infty$. Observe that $d(f,-\deg)=\deg(f)$. It is a fact that $f$ is proper if and only if $d(f,v)>0$ for all $v\in V_{\infty}.$

\smallskip

We now set\index{$f_* v$}
\begin{itemize}
\item
$f_* v :=0$ if $d(f,v) = 0$;
\item
$f_*v(P) =v(f^*P)$ if $d(f,v)> 0$.
\end{itemize}
In this way one obtains a valuation on $k[x,y]$ (that may be trivial); 
we then get a continuous map
$$f_{\d}:\{v\in V_{\infty}|\,\,d(f,v)>0\}\rightarrow V_{\infty}$$ defined by $$f_{\d}(v):=d(f,v)^{-1}f_*v ~.$$
This map extends to a continuous map $f_{\d}:\overline{\{v\in V_{\infty}|\,\,d(f,v)>0\}}\rightarrow V_{\infty}$. The image of any $v\in \partial{\{v\in V_{\infty}|\,\,d(f,v)>0\}}$ is a curve valuation defined by a rational curve with one place at infinity.

We now recall the following key result, \cite[Proposition 2.3\,,Theorem 2.4,\,Proposition 5.3.]{Favre2011}.
\begin{thm}\label{thm:eigenval}
There exists a valuation $v_*$
such that $\alpha(v_*)\geq 0\geq A(v_*)$, and
$f_*v_*=\la_1 v_*$.

If $\la_1(f)^2>\la_2(f)$, this valuation is unique .

If $\la_1(f)^2=\la_2(f)$, the set of such valuations is a closed segment in $V_{\infty}$.
\end{thm}

This valuation $v_*$ is called the {\em eigenvaluation} \index{eigenvaluation}of $f$ when $\la_1(f)^2>\la_2(f)$.

\section{The Green function of $f$}\label{sectiongreen}

\subsection{Subharmonic functions on $V_{\infty}$}
We refer to \cite[Section 3]{Xieb} for details.

To any $v\in V_\infty$ we attach its Green function\index{$g_v$}
$$
Z_v(w) := \alpha(v \wedge w)~.
$$
This is a decreasing continuous function taking values in $[-\infty, 1]$, satisfying
$g_v(-\deg) =1$.

Given any positive Radon measure $\rho$ on $V_\infty$ we define
$$
Z_\rho (w) := \int_{V_\infty} Z_v(w) \, d\rho(v)~.
$$
Observe that $g_v(w)$ is always well-defined as an element in $[-\infty, 1]$ since
$g_v \le 1$ for all $v$.

Then we recall the following result.
\begin{thm}[\cite{Xieb}]\label{thmrhotogrhoinj}
The map $\rho \mapsto Z_\rho$ is injective.
\end{thm}

One can thus make the following definition.
\begin{defi}
A function $\phi : V_\infty \to \R \cup \{-\infty\}$ is said to be \emph{subharmonic}\index{subharmonic} if there exists a positive Radon measure $\rho$ such that $\phi = Z_\rho$. In this case, we write
$\rho = \Delta \phi$ \index{$\Delta \phi$}and call it the \emph{Laplacian}\index{Laplacian} of $\phi$.
\end{defi}
%Denote by $\SH(V_{\infty})$\index{$\SH(V_{\infty})$} (resp. $\SH^{+}(V_{\infty})$)\index{$\SH^{+}(V_{\infty})$} the space of subharmonic functions on $V_\infty$ (resp. of non-negative subharmonic functions on $V_\infty$).

\bigskip

\subsection{Basic properties of the Green function of $f$}
We refer to \cite[Section 12]{Xiec} for details.

Let $f$ be a dominant polynomial endomorphism on $\A^2_k$ with $\la_1(f)^2>\la_2(f).$

By \cite[Section 12]{Xiec}, there exists a unique subharmonic function $\theta^*$ on $V_{\infty}$ such that
\begin{points}
\item $f^*\theta^*=\la_1\theta^*$;
\item $\theta^*(v)\geq 0$ for all $v\in V_{\infty}$;
\item $\theta^*(-\deg)=1$;
\item for all $v\in V_{\infty}$ satisfying $\alpha(v)>-\infty$, we have  $\theta^*(v)>0$ if and only if $d(f^n,v)>0$ for all $n\geq 0$ and $$\lim_{n\to \infty} f_{\d}^n(v)=v_*.$$
\end{points}

\section{Metrics on projective varieties defined over a valued field}\label{secmop}

A field with an absolute value is called a valued  field.
\medskip

\begin{defi}Let $(K,|\cdot|_v)$ be a valued field. For any integer $n\geq 1,$ we define a metric $d_v$ on the projective space $\mathbb{P}^n(K)$ by
$$d_v([x_0:\cdots:x_n],[y_0:\cdots:y_n])=\frac{\max_{0\leq i,j\leq n}|x_iy_j-x_jy_i|_v}{\max_{0\leq i\leq n}|x_i|_v\max_{0\leq j\leq n}|y_j|_v}$$ for any two points $[x_0:\cdots:x_n],[y_0:\cdots:y_n]\in \mathbb{P}^n(K).$
\end{defi}
\smallskip

Observe that when $|\cdot|_v$ is archimedean,  then the metric $d_v$ is not induced by a smooth riemannian metric. However it is equivalent to the restriction of the Fubini-Study metric on $\mathbb{P}^n(\mathbb{C})$ or $\mathbb{P}^n(\mathbb{R})$ to $\mathbb{P}^n(K)$ induced by any embedding $\sigma_v: K\hookrightarrow \R \text{ or } \C$.

\bigskip

More generally, for a projective variety $X$ defined over $K$, if we fix an embedding $\iota: X\hookrightarrow\mathbb{P}^n$, we may restrict the metric $d_v$ on $\mathbb{P}^n(K)$ to a metric $d_{v,\iota}$ on $X(K).$ This metric depends on the choice of embedding $\iota$  in general, but for different embeddings $\iota_1$ and $\iota_2$, the metrics $d_{v,\iota_1}$ and $d_{v,\iota_2}$ are equivalent.
Since we are mostly intersecting in the topology induced by these metrics, we shall usually write $d_v$ instead of $d_{v,\iota}$ for simplicity.

\newpage

\part{The existence of Zariski dense orbit}\label{part2}
The aim of this part is to prove Theorem \ref{thmconzhangatwo}.

\section{The attracting case}\label{sectionatt}
In this section, we prove Theorem \ref{thmconzhangatwo} in some special cases. In most of these cases, we find a Zariski dense orbit in some attracting locus. We also prove Theorem \ref{thmconzhangatwo} when $\la_1^2=\la_2>1$.

\medskip

Denote by $k$ an algebraically closed field of characteristic $0$. Let $f:\mathbb{A}^2\rightarrow \mathbb{A}^2$ be a dominant polynomial endomorphism defined over $k$. We have the following result.
\begin{lem}\label{lemlatglaodivzhang}If $\la_2^2(f)>\la_1(f)$ and the eigenvaluation $v_*$ is not divisorial,
then there exists a point $p\in \mathbb{A}^2(k)$ whose the orbit is Zariski dense in $\mathbb{A}^2$.
\end{lem}

\proof[Proof of Lemma \ref{lemlatglaodivzhang}] After replacing $k$ by an algebraically closed subfield of $k$ containing all the coefficients of $f$, we may suppose that the transcendence degree of $k$ over $\bar{\Q}$ is finite.

By \cite[Theorem 3.1]{Favre2011}, there exists a compactification $X$ of $\mathbb{A}^2$ defined over $k$ and a superattacting point $q\in X\setminus \A^2_k$ such that for any valuation $v\in V_{\infty}$ whose center in $X$ is $q$, we have $f_{\d}^nv\rightarrow v_*$ as $n\rightarrow \infty$.

By embedding $k$ in $\mathbb{C}$, we endow $X$ with the usual Euclidean topology. There exists a neighborhood $U$ of $q$ in $X$ such that
\begin{points}
\item
$I(f)\cap U=\emptyset$;
\item
$f(U)\subseteq U$;
\item
for any point $p\in U$, $f^n(p)\rightarrow q$ as $n\rightarrow\infty$.
\end{points}

Since $\mathbb{A}^2(k)$ is dense in $\mathbb{A}^2(\mathbb{C})$, there exists a point $p\in \mathbb{A}^2(k)\cap U.$
If the orbit of $p$ is not Zariski dense, then its Zariski closure $Z$ is a union of finitely many curves. Since $f^n(p)\rightarrow q$, $p$ is not preperiodic.
It follows that all the one dimensional irreducible components of $Z$ are periodic under $f$. Let $C$ be a one dimensional irreducible component of $Z$. Since there exists an infinite sequence $\{n_0<n_1<\cdots\}$ such that $f^{n_i}(p)\in C$ for $i\geq 0,$ we have $C$ contains $q=\lim_{i\rightarrow\infty}f^{n_i}(p)$. Then there exists a branch $C_1$ of $C$ at infinity satisfying $q\in C_1$. Thus $v_{C_1}$ is periodic which is a contradiction. It follows that $O(p)$ is Zariski dense.
\endproof

In many cases, for example $\la_1(f)^2>\la_2(f)$ and $v_*$ is divisorial,
there exists
a projective compactification $X$ of $\mathbb{A}^2$ and an irreducible component $E$ of $X\setminus \A^2$ satisfying $f_{\d}(v_E)=v_E$.
The following result proves Theorem \ref{thmconzhangatwo} when $f|_E$ is of infinite order.

\begin{lem}\label{lemconjzhangfenotid}Let $X$ be a projective compactification of $\mathbb{A}^2$ defined over $k$. Then $f$ extends to a rational selfmap on $X$. Let $E$ be an irreducible component of $X\setminus \A^2$ satisfying $f_{\d}(v_E)=v_E$, $d(f,v_E)\geq 2$ and $f^n|_E\neq \id$ for all $n\geq 1$. Then  there exists a point $p\in \mathbb{A}^2(k)$ whose orbit is Zariski dense in $\mathbb{A}^2$.
\end{lem}

\proof[Proof of Lemma \ref{lemconjzhangfenotid}]
There exists a finitely generated $\Z$-subalgebra $R$ of $k$ such that $X$, $E$ and $f$ are defined over the fraction field $K$ of $R$.

By \cite[Lemma 3.1]{Bell2006}, there exists a prime $\mathfrak{p}\geq 3$, an embedding of $K$ into $\Q_{\mathfrak{p}}$, and
a $\Z_{\mathfrak{p}}$-scheme $\mathcal{X}_{Z_{\mathfrak{p}}}$ such that the generic fiber is $X$, the specialization $E_{\mathfrak{p}}$ of $E$ is isomorphic to $\P^1_{\F_{\mathfrak{p}}}$, the specialization $f_{\mathfrak{p}}:X_{\mathfrak{p}}\dashrightarrow X_{\mathfrak{p}}$ of $f$ at the prime ideal ${\mathfrak{p}}$ of $Z_{\mathfrak{p}}$ is dominant, and $\deg f_{\mathfrak{p}}|_{E_{\mathfrak{p}}}=\deg f_E$.

\medskip

%We first suppose that $\deg (f_E)\geq 2.$

Since there are only finitely many points in the orbits of $I(f_{\mathfrak{p}})$ and the orbits of ramified points of $f_{\mathfrak{p}}$, by \cite[Proposition 5.5]{fa}, there exists a closed point $x\in E_{\mathfrak{p}}$ such that $x$ is periodic, $E_{\mathfrak{p}}$ is the unique irreducible component of $X_{\mathfrak{p}}\setminus \mathbb{A}^2_{\mathfrak{p}}$ containing $x$, $x\not\in I(f_{\mathfrak{p}}^n)$ for all $n\geq 0$ and $f_{\mathfrak{p}}|_{E_{\mathfrak{p}}}$ is not ramified at any point on the orbit of $x$. After replacing $\Q_{\mathfrak{p}}$ by a finite extension $K_{\mathfrak{p}}$ we may suppose that $x$ is defined over $O_{K_{\mathfrak{p}}}/\mathfrak{p}$.
After replacing $f$ by a positive iterate, we may suppose that $x$ is fixed by $f_{{\mathfrak{p}}}$.

The fixed point $x$ of $f_{\mathfrak{p}}$ defines an open and closed polydisc  $U$ in $X(K_{\mathfrak{p}})$ with respect to the ${\mathfrak{p}}$-adic norm $|\cdot|_{\mathfrak{p}}$. We have $f(U)\subseteq U$.
Observe that $f^*E= d(f,v)E$ in $U$ and $d(f,v)\geq 2$. So for all point $q\in U\cap \mathbb{A}^2(K_{\mathfrak{p}})$, we have $d_{\mathfrak{p}}(f^n(q),E)\rightarrow 0$ as $n\rightarrow \infty$.

Since $f_{\mathfrak{p}}|_{E_{\mathfrak{p}}}$ is not ramified at $x$, after replacing $f$ by some positive iterate, we may suppose that $df_{\mathfrak{p}}|_{E_{\mathfrak{p}}}(x)=1.$
By \cite[Theorem 1]{Poonen2014}, we have that for any point $q\in U\cap E$, there exists a ${\mathfrak{p}}$-adic analytic map $\Psi:O_{K_{\mathfrak{p}}} \rightarrow U\cap E$, such that for any $n\geq 0,$ we have $f^n(q)=\Psi(n).$

If there exists a preperiodic point $q$ in $U\cap E$, the there exists $m\geq 0$ such that $f^m(q)$ is periodic. Then there are infinitely many $n\in \mathbb{Z}^+\subseteq O_{K_{\mathfrak{p}}}$ such that $\Psi(n)=f^m(q)$. The fact that $O_{K_{\mathfrak{p}}}$ is compact shows that $\Psi$ is constant. It follows that $q$ is fixed. Thus all preperiodic points in $U\cap E$ are fixed by $f|_E$.

Let $S$ be the set of all fixed points in $U\cap E.$ Since $f^n|_E\neq \id$ for all $n\geq 1,$ $S$ is finite.

We first treat the case $S=\emptyset.$
Pick a point $p\in U\cap \mathbb{A}^2(\bar{\Q})$. Then $p$ is not preperiodic. If $O(p)$ is not Zariski dense, we denote by $Z$ its Zariski closure. Pick a one dimensional irreducible component $C$ of $Z$. We have $C\cap E\cap U\neq\emptyset$, and for all points $q\in C\cap E\cap U$, $q$ is preperiodic under $f|_E$. This contradicts our assumption, so $O(p)$ is Zariski dense.

Next we treat the case $S\neq \emptyset.$
Since $| df|_{E}(q_i)|_{\mathfrak{p}}=1$ for all $i=1,\cdots,m$, we have $df|_{E}(q_i)\neq 0$.
By embedding $K$ in $\C$, $X$ we endow $X$ with the usual Euclidean topology. By \cite[Theorem 3.1.4]{Abate2001}, for any $i=1,\cdots,m$, there exists a unique complex analytic manifold $W$ not contained in $E$  such that $f(W)=W$. It follows that
there are at most one irreducible algebraic curve $C_i\neq E$ in $X$ such that $q_i\in C_i$ and $f(C_i)\subseteq C_i$.
For convenience, if such an algebraic curve does not exists, we define $C_i$ to be $\emptyset$.

For any $n\geq 1$, by applying \cite[Theorem 3.1.4]{Abate2001} for $f^n$, if $C$ is a curve satisfying $q_i\in C$ and $f(C)\subseteq C$, then $C=C_i$. Moreover if $C'$ is an irreducible component of $f^{-1}(C_i)$ such that $q\in C'$, then for any point $y\in C'$ near $q$ w.r.t. the Euclidean topology, we have $f(p)\in C$. Then by \cite[(iv) of Theorem 3.1.4]{Abate2001}, we have $p\in C$. It follows that $C'=C$. Thus there exists a small open and closed neighborhood $U_i$ of $q_i$ w.r.t. to the norm $|\cdot|_{\mathfrak{p}}$ such that for all $j\neq i$, $q_j\not\in U_i$, $f(U_i)\subseteq U_i$ and $f^{-1}(C_i\cap U_i)\cap U_i=C_i\cap U_i$.

Observe that $\A^2(\bar{\Q}\cap K_{\mathfrak{p}})$ is dense in $\A^2(K_{\mathfrak{p}})$ w.r.t. $|\cdot|_{\mathfrak{p}}$.

There exists a $\bar{\Q}$-point $p$ in $U_1\setminus C_1(K_{\mathfrak{p}})$. If the orbit $O(p)$ of $p$ is not Zariski dense, denote by $Z$ its Zariski closure. Since $d_{\mathfrak{p}}(f^n(p),E)\rightarrow 0$ as $n\rightarrow \infty$, $p$ is not preperiodic. It follows that there exists a one dimensional irreducible component $C$ of $Z$ which is periodic. There exists $a,b>0$ such that $f^{an+b}(p)\in C$ for all $n\geq 0$.
Since $d_{\mathfrak{p}}(f^n(p),E)\rightarrow 0$ as $n\rightarrow \infty$ and $U_1$ is closed, there exists a point $q\in C\cap E \cap U_1$. It follows that $q$ is periodic under $f|_E$.Then $q$ is fixed and $q=q_1$. This implies $C=C_1$. Since $f^b(p)\in C_1\cap V$ and $f^{-1}(C_1)\cap V=C_1$, we have $p\in C_1$ which is a contradiction. It follows that $O(p)$ is Zariski dense.
%
%\medskip
%
%Then we treat the case $\deg f|_E=1$.
%
%By replacing $f$ by a suitable positive iteration, we may suppose that $f_{\mathfrak{p}}|_{E_{\mathfrak{p}}}=\id.$ There exists a closed point $x\in E_{\mathfrak{p}}$ such that $E_{\mathfrak{p}}$ is the unique irreducible component of $X_{\mathfrak{p}}\setminus \mathbb{A}^2_{\mathfrak{p}}$ containing $x$ and $x\not\in I(f_{\mathfrak{p}}^n)$ for all $n\geq 0$. By replacing $K$ by a finite extension we may suppose that $x$ is defined over $O_{K_{\mathfrak{p}}}/\mathfrak{p}$.
%The fix point $x$ of $f_{\mathfrak{p}}$ defines an open and closed polydisc $U$ in $X(K_{\mathfrak{p}})$ w.r.t. the norm $|\cdot|_{\mathfrak{p}}$ such that $f(U)\subseteq U$.
%Since $f|_E\neq \id$, there are at most two preperiodic point under $f_E$. So we may suppose that there are no preperiodic points in $U\cap E$.
%Observe that $f^*E= d(f,v)E$ in $U$ and $d(f,v)\geq 2$. For all point $p\in U\cap \mathbb{A}^2(K_{\mathfrak{p}})$, we have $d_{\mathfrak{p}}(f^n(p),E)\rightarrow 0$ as $n\rightarrow \infty.$ Pick a point $p\in U\cap \mathbb{A}^2(\bar{\Q})$, then $p$ is not preperiodic. If $O(p)$ is not Zariski dense, we denote by $Z$ the Zariski closure of $p$. Pick $C$ a one dimensional irreducible component of $Z$, we have $C\cap E\cap U\neq\emptyset$ and for all point $q\in C\cap E\cap U$, $q$ is preperiodic under $f|_E$. It contradicts our assumption, so $O(p)$ is Zariski dense.
\endproof

\begin{pro}\label{proconjofzhanglatelaosgeqt}If $\la_2^2(f)=\la_1(f)>1$ then Theorem \ref{thmconzhangatwo} holds.
\end{pro}
\proof[Proof of Proposition \ref{proconjofzhanglatelaosgeqt}]
By \cite[Proposition 5.1]{Favre2011} and \cite[Proposition 5.3]{Favre2011} there exists a divisorial valuation  $v_*\in V_{\infty}$, satisfying $f_{\d}(v_*)=v_*$ and $d(f,v_*)=\la_1\geq 2.$ Moreover, there exists a compactification $X$ of $\A^2$ and an irreducible component $E$ in $X\setminus \A^2$ satisfying $v_*=v_E$ and $\deg (f|_E)=\la_1\geq 2$. We conclude our theorem by invoking Lemma \ref{lemconjzhangfenotid}.
\endproof

\section{Totally invariant curves}\label{sectiontot}
Let $f:\mathbb{A}^2\rightarrow \mathbb{A}^2$ be a dominant polynomial endomorphism defined over an algebraically closed field $k$ of of characteristic $0$.
Let $X$ be a compactification of $\mathbb{A}^2_{k}$. Then $f$ extends to a rational selfmap on $X$.

As in \cite{Invarianthypersurfaces}, a curve $C$ in $X$ is said to be $totally$ $invariant$ if the strict transform $f^{\#}C$ equals $C$.

If there are infinitely many irreducible totally invariant curves in $\mathbb{A}^2_k$, then \cite[Theorem B]{Invarianthypersurfaces} shows that $f$ preserves a nontrivial fibration:
\begin{pro}\label{prototinvacur}If there are infinitely many  irreducible curves in $\mathbb{A}^2_{k}$ that are totally invariant under $f$, then there is a nonconstant
rational function $g$ satisfying $g\circ f= g$.
\end{pro}
In this section, we give a direct proof of this result.

\proof[Proof of Proposition \ref{prototinvacur}]Let $\{C_i\}_{i\geq 1}$ be an infinite sequence of distinct irreducible totally invariant curves in $\mathbb{A}^2_k$. Since the ramification locus of $f$ is of dimension at most one,
 after replacing $\{C_i\}_{i\geq 1}$ by an infinite subsequence, we may suppose that $C_i$ is not contained in the ramification locus of $f$ for any $i\geq 0$. Then we have $\ord_{C_i}f^*C_i=1$ for all $i\geq 0$.

Let $E_1,\cdots,E_s$ be the set of irreducible curves in $\mathbb{A}^2$ contracted by $f$. Let $V$ be the $\Q$-subspace in $\Div(\A^2)\otimes \Q$ spanned by $E_i$, $i=1,\cdots s.$ Then we have $f^*C_i=C_i+F_i$ where $F_i\in V$ for all $i\geq 1.$ Set $W_i:=\bigcap_{j\geq i}(\sum_{t\geq j}\mathbb{Q}F_t)\subseteq V$. We have $W_{i+1}\subseteq W_i$ for $i\geq 1.$ Since $V$ is of finite dimension, there exists $n_0\geq 1$, such that $W_i=W_{n_0}$ for all $i\geq n_0$. We may suppose that $n_0=1$ and set $W:=W_{n_0}$. Moreover we may suppose that $W$ is generated by $F_1,\cdots,F_l$ where $l=\dim W$.

For all $i\geq l+1$, we have $F_i=\sum_{j=1}^la^i_jF_j$ where $a^i_j\in \mathbb{Q}.$ Since $F_i=f^*C_i-rC_i$, we have $f^*(C_i-\sum_{j=1}^la^i_jC_j)=r(C_i-\sum_{j=1}^la^i_jC_j).$ There exists $n_i\in \mathbb{Z}^+$ such that $n_ia^i_j\in \mathbb{Z}$ for all $j=1,\cdots,l$. Then we have $f^*(n_iC_i-\sum_{j=1}^ln_ia^i_jC_j)=r(n_iC_i-\sum_{j=1}^ln_ia^i_jC_j).$
Up to multiplication by a nonzero constant, there exists a unique $g_i\in k(x,y)\setminus \{0\}$ such that $\Div(g_i)=n_iC_i-\sum_{j=1}^ln_ia^i_jC_j.$
It follows that
$f^*g_i=A_ig_i$ where $A_i\in k\setminus \{0\}.$ Since $n_iC_i-\sum_{j=1}^ln_ia^i_jC_j\neq 0$ for $i\geq l+1$, $g_i$ is non constant for $i\geq l+1.$
This concludes the proof.
\endproof

\begin{cor}\label{corfbirconz}If $f$ is birational, then either there is a nonconstant
rational function $g$ satisfying $g\circ f= g$,
or there exists a point $p\in \mathbb{A}^2(k)$ with Zariski dense orbit.
\end{cor}
\proof[Proof of Corollary \ref{corfbirconz}]
There exists a finite generated $\Q$-subalgebra $R$ of $k$ such that $f$ is defined over $R$. Denote by $K$ the fraction field of $R$.

By \cite[Lemma 3.1]{Bell2006}, there exists a prime $\mathfrak{p}\geq 3$, and an embedding of $R$ into $\Z_{\mathfrak{p}}$ such that all coefficients of $f$ are of $\mathfrak{p}$-adic norm $1$. Denote by $\mathbb{F}$ the algebraic closure of $\F_\mathfrak{p}$. Then the degree of the specialization $f_{\mathfrak{p}}:\mathbb{A}^2_{\mathbb{F}}\dashrightarrow \mathbb{A}^2_{\mathbb{F}}$ of $f$ equals $\deg f$.

By \cite[Proposition 6.2]{Xie2015}, there exists a noncritical periodic point $x\in \mathbb{A}^2(\mathbb{F})$. After replacing $\Q_p$ by a finite extension $K_{\mathfrak{p}}$, we may suppose that $x$ is defined over $O_{K_\mathfrak{p}}/\mathfrak{p}$. Replacing  $f$ by a suitable iterate, we may suppose that $x$ is fixed. Since $x$ is noncritical, $df_{\mathfrak{p}}(x)$ is invertible. After replacing $f$ by a suitable iterate, we may suppose that $df_{\mathfrak{p}}(x)=\id$.

The fixed point $x$ defines an open and closed neighborhood   $U$ in $\mathbb{A}^2(K_{\mathfrak{p}})$ with respect to $d_{\mathfrak{p}}$ such that $f(U)\subseteq U$.
By applying \cite[Theorem 1]{Poonen2014}, we have that for any point $q\in U$, there exists a ${\mathfrak{p}}$-adic analytic map $\Psi:O_{K_{\mathfrak{p}}} \rightarrow U$ such that
for any $n\geq 0,$ we have $f^n(q)=\Psi(n).$

Arguing by contradiction, we suppose that the orbit $O(q)$ of $q$ is not Zariski dense for any $q\in \A^2(k\cap K_{\mathfrak{p}})$. As in the proof of Lemma \ref{lemconjzhangfenotid}, if $q$ is preperiodic, then $q$ is fixed. Suppose that $q$ is non-preperiodic. There exists an irreducible curve $C$ such that $f^n(q)\in C$ for infinitely many $n\geq 0.$ Let $P$ be a polynomial such that $C$ is defined by $P=0$. Then $P\circ \Psi$ is an analytic function on $O_{K_{\mathfrak{p}}}$ having infinitely many zeros. It follows that $P\circ \Psi\equiv 0$ and then $f^n(q)\in C$ for all $n\geq 0$. Then we have $f(C)=C$.
Since $f$ is birational, a curve $C$ is totally invariant by $f$ if and only if $f(C)=C$. We may suppose that $f\neq \id$. Since $\overline{\mathbb{Q}}\cap K_{\mathfrak{p}}\subseteq k$ and the $\overline{\mathbb{Q}}\cap K_{\mathfrak{p}}$- points in $U$ are Zariski dense in $\mathbb{A}^2_{K_{\mathfrak{p}}}$, there are infinitely many irreducible totally invariant curves in $\mathbb{A}^2$. We conclude our Corollary by invoking Proposition \ref{prototinvacur}.
\endproof

\section{Proof of Theorem \ref{thmconzhangatwo}}\label{sectionproof}
Let $f:\mathbb{A}^2\rightarrow \mathbb{A}^2$ be a dominant polynomial endomorphism defined over an algebraically closed field $k$ of characteristic $0$.
%\proof[Proof of Theorem \ref{thmconzhangatwo}]

After replacing $k$ by an algebraically closed subfield which contains all the coefficients of $f$, we may suppose that the transcendence degree of $k$ over $\bar{\Q}$ is finite.

By Lemma \ref{lemlatglaodivzhang}, Proposition \ref{proconjofzhanglatelaosgeqt} and Corollary \ref{corfbirconz}, we may suppose that
$\la_1^2>\la_2>1$ and $v_*$ is divisorial. Suppose that $v_*=v_E$ for some irreducible exceptional divisor $E$ in some compactification $X$. If $f^n|_E\neq \id$ for all $n\geq 1$, Lemma \ref{lemconjzhangfenotid} concludes the proof. So after replacing $f$ by a suitable iterate, we may suppose that $f|_E=\id.$

By choosing a suitable compactification $X\in \mathcal{C}_0$, we may suppose that $E\cap I(f)=\emptyset$.
There exists a subfield $K$ of $k$ which is finite generated over $\Q$ such that $X, f, E, I(f)$ are defined over $K$. Moreover we may suppose that $E\simeq \P^1$ over $K.$

By \cite[Lemma 3.1]{Bell2006}, there exists a prime $\mathfrak{p}\geq 3$, such that we can embed $K$ into $\Q_{\mathfrak{p}}$.
Further there exists an open and closed set $U$ of $X(\Q_{\mathfrak{p}})$ w.r.t. the norm $|\cdot|_{\mathfrak{p}}$ containing $E$ and satisfying $U\cap I(f)=\emptyset$, $f(U)\subseteq U$ and $d_{\mathfrak{p}}(f^{n}(p),E)\rightarrow 0$ as $n\rightarrow \infty$ for any $p\in U$.

By \cite[Theorem 3.1.4]{Abate2001}, for any point $q\in E(\bar{K}\cap Q_{\mathfrak{p}})$, there is at most one irreducible algebraic curve $C_q\neq E$ in $X$ such that $q\in C_q$ and $f(C_q)\subseteq C_q$.
For convenience, if such an algebraic curve does not exists, set $C_q:=\emptyset$. Further if $C_q\neq \emptyset$, $C_q$ is smooth at $q$ and intersects $E$ transitively.

If there are only finitely many points $q\in E(\bar{K}\cap \Q_{\mathfrak{p}})$ such that $C_q$ is an algebraic curve, there exists a point $q\in E$ and an open and closed set $V$ w.r.t. the norm $|\cdot |_{\mathfrak{p}}$ containing $q$ such that $f(V)\subseteq V$ and $C_x=\emptyset$ for all $x\in V\cap E(\bar{K}\cap \Q_{\mathfrak{p}})$. Pick a point $p\in V\cap \mathbb{A}^2(\bar{K}\cap Q_{\mathfrak{p}})$, then $p$ is not preperiodic. If $O(p)$ is not Zariski dense, we denote by $Z$ its Zariski closure. There exists a one dimensional irreducible component $C$ of $Z$ which is periodic under $f$. Since $C\cap O(p)$ is infinite, $C$ is definied over $\bar{K}\cap Q_{\mathfrak{p}}$. Since $d_v(f^{n}(p),E)\rightarrow 0$ as $n\rightarrow \infty$, we have $C\cap E(\bar{K}\cap Q_{\mathfrak{p}})\cap V\neq\emptyset$. Pick $q\in C\cap E(\bar{K}\cap Q_{\mathfrak{p}})\cap V$, then we have that $C_q=C$ is an algebraic curve which contradicts our assumption. Thus $O(p)$ is Zariski dense.

Otherwise there exists an infinite sequence of points $q_i$, $i\geq 1$ such that $C_i:=C_{q_i}$ is an algebraic curve.
%
%
%Recall Siegel's Theorem (see \cite{Hindry2000} for details)
%
%\begin{thm}[Siegel's Theorem]\label{thmsiegel}
%Let $C$ be a curve over a number field $K$ and $g\in K(C)$ be
%a nonconstant rational function on $C$. If either $C$ is not rational or
%$g$ has at least three distinct poles, then the set $\{p\in C(K)|g(p)\in O_{K,S}\}$
%is finite.
%\end{thm}
%
%
Since $f|_{C_i}$ is an endomorphism of $C_i$ of degree $\la_1>1$, every $C_i$ is rational and has at most two branches at infinity. We may suppose that $E$ is the unique irreducible component of $X\setminus \mathbb{A}^2$ containing $q_i$ for all $i\geq 1$ and $C_i\neq C_j$ for $i\neq j$. We need the following result, which is proved below.
\begin{lem}\label{lemconjzhangdegbound}After replacing $\{C_i\}_{i\geq 1}$ by an infinite subsequence, we have that either $\deg (C_i)$ is bounded or Theorem \ref{thmconzhangatwo} holds.
\end{lem}
Suppose that $\deg (C_i)$ is bounded. Pick an ample line bundle $L$ on $X$. Then there exists $M>0$ such that $(C_i\cdot L)\leq M$ for all $i\geq 1$.

There exist a smooth projective surface $\Gamma$, a birational morphism $\pi_1:\Gamma\rightarrow X$ and morphism $\pi_2:\Gamma\rightarrow X$ satisfying $f=\pi_2\circ\pi_1.$ We denote by $f_*$ the map $\pi_{2*}\circ\pi^*_1:\Div X\rightarrow \Div X$. Let $E_{\pi_1}$ be the union of exceptional irreducible divisors of $\pi_1$ and $\mathfrak{E}$ be the set of effective divisors in $X$ supported by $\pi_2(E_{\pi_1}).$ It follows that for any curve $C$ in $X$, there exists $D\in \mathfrak{E}$ such that $f_*C=\deg(f|_C)f(C)+D$.

For any effective line bundle $M\in \Pic(X)$, the projective space $H_M:=\mathbb{P}(H^0(M))$ parameterizes the curves $C$ in the linear system $|M|$.
Since $\Pic^0(X)=0$, for any $l\geq 0$, there are only finitely many effective line bundles satisfying $(M\cdot L)\leq l.$

Then $H^l:=\coprod_{(M\cdot L)\leq l}H_M$ is a finite union of projective spaces and
it parameterizes the curves $C$ in $X$ satisfying $(C\cdot L)\leq l$.

There exists $d\geq 1$ such that $dL-f^*L$ is nef. Then, for any curve $C$ in $X$, we have $(f_*C\cdot L)=(C\cdot f^*L)\leq d(C\cdot L).$
It follows that $f_*$ induces a morphism $F:H^l\rightarrow H^{dl}$ by $C\rightarrow f_*C$ for all $l\geq 1.$ For all $l\geq 1$, $a\in \mathbb{Z}^+$ and $D\in \mathfrak{E}$, there exists an embedding $i_{a,D}:H_l\rightarrow H_{al+(D\cdot L)}$ by $C\mapsto aC+D$. Let $Z_1,\cdots,Z_m$ be all irreducible components of the Zariski closure of $\{C^j\}_{j\leq -1}$ in $H^M$ of maximal dimension. For any $i\in\{1,\cdots,m\}$, there exists $l\leq M$ such that $(C\cdot L)=l$ for all $C\in Z_i$. Let $S$ be the finite set of pairs $(a,D)$ where $a\in \mathbb{Z}^+$, $D\in \mathfrak{E}$ satisfying $al+(D\cdot L)\leq dM$. Then we have $F(Z_i)\subseteq \bigcup_{j=1,\cdots,m}\bigcup_{(a,D)\in S}i_{a,D}(Z_j)$. It follows that there exists a unique $j_i\in \{1,\cdots,m\}$, and a unique $(a,D)\in S$ such that $F(Z_i)=i_{a,D}(Z_{j_i}).$ Observe that, the map $i\mapsto j_i$ is an one to one map of $\{1,\cdots,m\}$. After replacing $f$ by a positive iterate, we may suppose that $j_i=i$ and $F(Z_i)\subseteq i_{a_{Z_i},D_{Z_i}}Z_i$ for all $i=1,\cdots,m.$ Set $Z:=Z_1$, $a=a_{Z_1}$, $D=D_{Z_1}$. We may suppose that $C_i\in Z$ for all $i\geq 1.$
Since $f(C_i)=C_i$, for all $i\geq 1$, we have $i_{a_{Z_1},D_{Z_1}}^{-1}\circ F|_Y=\id$.

For any point $t\in Z$, denote by $C^t$ the curve parameterized by $t$. Then $f(C^t)= C^t$ for all $t\in Z$. Since $C_i$ intersects $E$ transversely at at most two points, there exists $s\in \{1,2\}$ such that $C^t$ intersects $E$ transversely at $s$ points for a general $t\in Z$.
It follows that, for a general point $t\in Z$, there are exactly $s$ points $q^1_t,\cdots,q_t^s \in E$ such that $C^t=C_{q_t^j}$ for $j=1,\cdots,s$.

Set $Y:=\{(p,t)\in X\times Z|\,\,p\in C^t\}.$ Denote by $\pi_1:Y\rightarrow X$ the projection to the first coordinate and by $\pi_2: Y\rightarrow Z$ the projection to the second coordinate. Since $C^t\cap E$ is not empty for general $t\in Z$, the map $\pi_2|_{\pi_1^{*}E}$ is dominant. We see that $f$ induces a map $T:Y\rightarrow Y$ defined by $(p,t)\rightarrow (f(p),t)$.
Since there are infinitely many points in $E$ contained in $\pi_1(\pi_1^{*}E)$, so $\pi_1|_{\pi_1^{*}E}:\pi_1^{*}E\rightarrow E$ is dominant. For a general point $t\in Z$, there are exactly $s$ points in $\pi_1^{*}E$. It follows that the map $\pi_2|_{\pi_1^{*}E}$ is generically finite of degree $s$. For a general point $q\in E$, there exists only one point $(q,C_q)\in \pi_1^{*}E$. Hence $\pi_1(\pi_1^{*}E):\pi_1^{*}E\rightarrow E$ is birational. It follows that $Z$ is a rational curve and $Y$ is a surface. Thus the morphism $\pi_1:Y\rightarrow X$ is generically finite. Let $p$ be a general point in $\mathbb{A}^2$. If $\#\pi_1^{-1}(p)\geq 2$, then there are $t_1\neq t_2\in Z$ such that $p\in C^{t_1}\cap C^{t_2}.$ Since there exists $M'>0$ such that $\deg C^{t}\leq M'$ for all $t\in Z$, we have $\#(C_{t_1}\cap C_{t_2})\leq M'^2$. If follows that there exist $a<b\in \{0,\cdots, M'^2\}$ such that $f^a(p)=f^b(p)$. This contradicts the assumption that $p$ is general. It follows that the morphism $\pi_1:Y\rightarrow X$ is birational. Identify $Z$ with $\mathbb{P}^1$.
Set $g:=\pi_2\circ \pi_1^{-1}$. Then we have $g\circ f=g$, which completes the proof.
\endproof

\proof[Proof of Lemma \ref{lemconjzhangdegbound}]
If $C_i$ has only one place at infinity for all $i\geq 1$, then $\deg C_i=b_E$. So we may suppose that $C_i$ has two places at infinity for all $i\geq 1$.

\smallskip

We first suppose that $\#C_i\cap E=2$ for all $i\geq 1.$ We may suppose that $C_i\cap E\cap \Sing (X\setminus \mathbb{A}^2)=\emptyset$ for all $i\geq 1.$ Then for all $i\geq 1$, we have $\deg C_i=2b_E$.

\smallskip

Then we may suppose that $C_i\cap E=\{q_i\}$ for all $i\geq 1.$ Let $c_i$ be the unique branch of $C$ at infinity centered at $q_i$ and $w_i$ the unique branch of $C$ at infinity not centred at $q_i$. Since $f(C_i)=C_i$, we have $f_{\d}(c_i)=c_i$ and $f_{\d}(w_i)=w_i$.

\smallskip

We first treat the case $\theta^*=Z_{v^*}$ where $v^*\in V_{\infty}$ is divisorial. It follows that $\la_2/\la_1\in \mathbb{Z}^+$.
Observe that $\deg (f|_{C_i})=\la_1$ for $i$ large enough.

If $\la_1=\la_2$, then the strict transform $f^{\#}(C_i)$ equals to $C_i$ for $i$ large enough, and the lemma follows from Proposition \ref{prototinvacur}.

Otherwise we have $\la_2/\la_1\geq 2$. Set $v^*:=v_{E'}$. Then we have $d(f,v^*)=\la_2/\la_1\geq 2$ and $\deg (f|_{E'})=\la_1>1$.  By Lemma \ref{lemconjzhangfenotid}, Theorem \ref{thmconzhangatwo} holds.

\smallskip

Next we treat the case $\theta^*=Z_{v^*}$ where $v^*\in V_{\infty}$ is not divisorial. Since $\alpha(v^*)=0$, $v^*$ can not be irrational. Thus $v^*$ is infinitely singular, and hence an end in $V_{\infty}.$ Then $f$ is proper.

By \cite[Proposition 15.2]{Xiec}, there exists $v_1<v^*$ such that for any valuation $v\neq v^*$ in $U:=\{w\in V_{\infty}|\,\,w>v_1\}$, there exists $N\geq 1$ such that $f_{\d}^n(v)\not\in U$ for all $n\geq N.$
It follows that there is no curve valuation in $U$ which is periodic under $f_{\d}$. Let $U_E$ be the open set in $V_{\infty}$ consisting of all valuations whose center in $X$ is contained $E$. If follows that $w_i\not\in U\cup U_E$ for all $i\geq 1$. Hence $w_i\not\in U\cup f_{\d}^{-N}(U_E)$ for all $N\geq 0.$ Set $W_{-4}:=\{v\in V_{\infty}|\,\, \alpha(v)\geq -4\}$. By \cite[Proposition 11.6]{Xiec} and the fact that $W_{-4}\setminus U$ is compact, there exists $N\geq 0$ such that $f_{\d}^N(W_{-4}\setminus U)\subseteq U_E$. Then we have $W_{-4}\subseteq U \cup f^{-N}_{\d}(U_E).$

Since the boundary $\partial (V_{\infty}\setminus (U \cup f^{-N}_{\d}(U_E)))$ of $V_{\infty}\setminus (U \cup f^{-N}_{\d}(U_E))$ is finite and $w_i\in V_{\infty}\setminus (U \cup f^{-N}_{\d}(U_E))$ for all $i\geq 1$, we may suppose that there exists $w\in \partial (V_{\infty}\setminus (U \cup f^{-N}_{\d}(U_E)))$ satisfying $w_i>w$ for all $i\geq 1.$

If, for all $i\geq 1$, we have $(w_i\cdot l_{\infty})\leq 1/2\deg (C_i)$, then $\deg C_i=(w_i\cdot l_{\infty})+ (c_i\cdot l_{\infty})\leq 1/2\deg(C_i)+b_E$. It follows that $\deg(C_i)\leq 2b_E$.

Thus we may suppose that $(w_i\cdot l_{\infty})\geq 1/2\deg (C_i)$ for all $i\geq 1.$
For any $i\neq j$, the intersection number $(C_i\cdot C_j)$ is the sum of the local intersection numbers at all points in $C_i\cap C_j$. Since all the local intersection numbers are positive, we have
$$\deg(C_i)\deg(C_j)\geq(c_i\cdot c_j)+(w_i\cdot w_j).$$
By the calculation in Section \ref{sectionlocalinter}, we have
$$(c_i\cdot c_j)+(w_i\cdot w_j)$$$$\geq b_E^2(1-\alpha(v_E))+(w_i\cdot l_{\infty})(w_j\cdot l_{\infty})(1-\alpha(w))$$$$\geq b_E^2(1-\alpha(v_E))+5/4\deg(C_i)\deg(C_j).$$ Thus $\deg(C_i)\deg(C_j)\leq -4b_E^2(1-\alpha(v_E))<0$, which is a contradiction.

\smallskip

Finally, we treat the case when $\#\Supp\Delta\theta^*\geq 2$. By \cite[Theorem 2.4]{Favre2011}, we have $\theta^*>0$ on  the set $W_0:=\{v\in V_{\infty}|\,\, \alpha(v)\geq 0\}$.
Set $W_{-1}:=\{v\in V_{\infty}|\,\, \alpha(v)\geq -1\}$ and $Y:=\{v\in W_{-1}|\,\, \theta^*(v)=0\}$. By \cite[Proposition 11.2]{Xiec}, $Y$ is compact. For any point $y\in Y$, there exists $w_y< y$ satisfying $\alpha(w_y)\in (-1,0)$. Set $U_y:=\{v\in V_{\infty}|\,\, v>w_y\}$. There are finitely many points $y_1,\cdots,y_l$ such that $Y\subseteq \bigcup_{i=1}^lU_{y_i}$. Pick $r:=1/2\min\{-\alpha(w_{y_i})\}_{i=1,\cdots,l}$. Then $r\in (0,1)$, and $W_{-r}\cap (\bigcup_{i=1}^lU_{y_i})=\emptyset$. It follows that there exists $t>0$ such that $\theta^*\geq t$ on $W_{-r}.$ By \cite[Proposition 11.6]{Xiec}, there exists $N\geq 0$ such that $f^N_{\d}(W_{-r})\subseteq U_E$.
Then we have $W_{-r}\subseteq f^{-N}_{\d}(U_E).$

Since the boundary $\partial (V_{\infty}\setminus (U \cup f^{-N}_{\d}(U_E)))$ of $V_{\infty}\setminus (U \cup f^{-N}_{\d}(U_E))$ is finite and $w_i\in V_{\infty}\setminus (U \cup f^{-N}_{\d}(U_E))$ for all $i\geq 1$, we may suppose that there exists $w\in \partial (V_{\infty}\setminus (U \cup f^{-N}_{\d}(U_E)))$ satisfying $w_i>w$ for all $i\geq 1.$

Pick $\delta\in (0, \frac{r}{2(1+r)}).$
If, for all $i\geq 1$, we have $(w_i\cdot l_{\infty})\leq (1-\delta)\deg (C_i)$, then $\deg C_i=(w_i\cdot l_{\infty})+ (c_i\cdot l_{\infty})\leq (1-\delta)\deg(C_i)+b_E$. It follows that $\deg(C_i)\leq b_E/\delta$.

Thus we may suppose that $(w_i\cdot l_{\infty})\geq (1-\delta)\deg (C_i)$ for all $i\geq 1.$
For any $i\neq j$, we have $$\deg(C_i)\deg(C_j)\geq(c_i\cdot c_j)+(w_i\cdot w_j)$$$$\geq b_E^2(1-\alpha(v_E))+(w_i\cdot l_{\infty})(w_j\cdot l_{\infty})(1-\alpha(w))$$$$\geq b_E^2(1-\alpha(v_E))+(1-\delta)^2(1+r)\deg(C_i)\deg(C_j).$$ Set $t:=(1-\delta)^2(1+r)-1$, then we have
$$t>(1-2\delta)(1+r)-1>(1-r/(1+r))(1+r)-1=0,$$
and hence $\deg(C_i)\deg(C_j)\leq -t^{-1}b_E^2(1-\alpha(v_E))<0$ which is a contradiction. This concludes the proof.
\endproof

\bibliography{dd}

\begin{thebibliography}{10}

\bibitem{Abate2001}
{M}arco {A}bate.
\newblock {\em {A}n {I}ntroduction to {H}yperbolic {D}ynamical {S}ystems}.
\newblock {I}stituti {E}ditoriali e {P}oligrafici {I}nternazionali, 2001.

\bibitem{Amerik}
Ekaterina Amerik.
\newblock Existence of non-preperiodic algebraic points for a rational self-map
  of infinite order.
\newblock {\em Math. Res. Lett.}, 18(2):251--256, 2011.

\bibitem{E.Amerik2011}
Ekaterina Amerik, Fedor Bogomolov, and Marat Rovinsky.
\newblock Remarks on endomorphisms and rational points.
\newblock {\em Compositio Math.}, 147:1819--1842, 2011.

\bibitem{Amerik2008}
Ekaterina Amerik and Fr{\'e}d{\'e}ric Campana.
\newblock Fibrations m\'eromorphes sur certaines vari\'et\'es \`a fibr\'e
  canonique trivial.
\newblock {\em Pure Appl. Math. Q.}, 4(2, part 1):509--545, 2008.

\bibitem{Bell2006}
Jason~P. Bell.
\newblock A generalised {S}kolem-{M}ahler-{L}ech theorem for affine varieties.
\newblock {\em J. London Math. Soc. (2)}, 73(2):367--379, 2006.

\bibitem{Bell}
Jason~Pierre Bell, Dragos Ghioca, and Thomas~John Tucker.
\newblock Applications of {$p$}-adic analysis for bounding periods for
  subvarieties under \'etale maps.
\newblock {\em Int. Math. Res. Not. IMRN}, (11):3576--3597, 2015.

\bibitem{Invarianthypersurfaces}
Serge Cantat.
\newblock Invariant hypersurfaces in holomorphic dynamics.
\newblock {\em Math. Research Letters}, 17(5):833--841, 2010.

\bibitem{fa}
Najmuddin Fakhruddin.
\newblock Questions on self maps of algebraic varieties.
\newblock {\em J. Ramanujan Math. Soc.}, 18(2):109--122, 2003.

\bibitem{Fakhruddin2014}
Najmuddin Fakhruddin.
\newblock The algebraic dynamics of generic endomorphisms of {$\Bbb{P}^n$}.
\newblock {\em Algebra Number Theory}, 8(3):587--608, 2014.

\bibitem{Favre2004}
Charles Favre and Mattias Jonsson.
\newblock {\em The valuative tree}, volume 1853 of {\em Lecture Notes in
  Mathematics}.
\newblock Springer-Verlag, Berlin, 2004.

\bibitem{Favre2007}
Charles Favre and Mattias Jonsson.
\newblock Eigenvaluations.
\newblock {\em Ann. Sci. \'Ecole Norm. Sup. (4)}, 40(2):309--349, 2007.

\bibitem{Favre2011}
Charles Favre and Mattias Jonsson.
\newblock Dynamical compactifications of $\bold {C}^2$.
\newblock {\em Ann. of Math. (2)}, 173(1):211--248, 2011.

\bibitem{Jonsson}
Mattias Jonsson.
\newblock {D}ynamics on {B}erkovich spaces in low dimensions.
\newblock In J.~Nicaise A.~Ducros, C.~Favre, editor, {\em {B}erkovich spaces
  and applications}, volume 2119 of {\em Lecture Notes in Mathematics}.
  Springer, 2015.

\bibitem{Medvdevv1}
A.~Medvedev and T.~Scanlon.
\newblock Polynomial dynamics.
\newblock arXiv:0901.2352v1.

\bibitem{Medvdev}
Alice Medvedev and Thomas Scanlon.
\newblock Invariant varieties for polynomial dynamical systems.
\newblock {\em Ann. of Math. (2)}, 179(1):81--177, 2014.

\bibitem{Poonen2014}
Bjorn Poonen.
\newblock {$p$}-adic interpolation of iterates.
\newblock {\em Bull. Lond. Math. Soc.}, 46(3):525--527, 2014.

\bibitem{Xiec}
Junyi Xie.
\newblock The {D}ynamical {M}ordell-{L}ang {C}onjecture for polynomial
  endomrophisms on the affine plane.
\newblock arXiv:1503.00773.

\bibitem{Xieb}
Junyi Xie.
\newblock When the intersection of valuation rings of $k[x,y]$ has
  transcendence degree 2?
\newblock arXiv:1403.6052.

\bibitem{Xie2015}
Junyi Xie.
\newblock Periodic points of birational transformations on projective surfaces.
\newblock {\em Duke Math. J.}, 164(5):903--932, 2015.

\bibitem{zhang}
Shou-Wu Zhang.
\newblock {D}istributions in {A}lgebraic {D}ynamics.
\newblock {\em {S}urveys in differential geometry}, 10:381--430, International
  Press 2006.

\end{thebibliography}

\end{document}